\documentclass[11pt]{paper}

\usepackage{amsmath} %[intlimits]
\usepackage{amssymb}
\usepackage{amsfonts}
\usepackage{latexsym}
\usepackage{amsthm}
\usepackage{bbm}
\usepackage{mathrsfs}
\usepackage{bm}
\usepackage{xcolor}
\usepackage{subcaption}
\usepackage{graphicx, color}

\definecolor{darkred}{rgb}{0.9,0,0.3}
\definecolor{darkblue}{rgb}{0,0.3,0.9}

\theoremstyle{plain} %plain, definition, remark
\newtheorem{theorem}{Theorem}[section]
\newtheorem*{theorem*}{Theorem}

\newtheorem*{lemma*}{Lemma}

\newtheorem*{corollary*}{Corollary}
\newtheorem{proposition}[theorem]{Proposition}
\newtheorem*{proposition*}{Proposition}
\newtheorem{definition}[theorem]{Definition}
\newtheorem*{definition*}{Definition}

\newtheorem*{conjecture*}{Conjecture}

\theoremstyle{definition} %plain, definition, remark
\newtheorem{remark}[theorem]{Remark}
\newtheorem*{remark*}{Remark}

\renewcommand{\H}{\mathbb{H}}

\renewcommand{\leq}{\leqslant}
\renewcommand{\geq}{\geqslant}
\renewcommand{\epsilon}{\varepsilon}

\begin{document}
\title{A new approach to SLE phase transition}

\author{Dmitry Beliaev, Terry J. Lyons, Vlad Margarint\thanks{ TJL and VM would like to acknowledge the support of ERC (Grant Agreement No.291244 Esig), DB and VM were partially funded by EPSRC Grant EP/M002896/1}}

\maketitle

\begin{abstract}
It is well know that $SLE_\kappa$ curves exhibit a phase transition at $\kappa=4$. For $\kappa\le 4$ they are simple curves with probability one, for $\kappa>4$ they are not. The standard proof is based on the analysis of the Bessel SDE of dimension $d=1+4/\kappa$. We propose a different approach which is based on the analysis of the Bessel SDE with $d=1-4/\kappa$. This not only gives a new perspective, but also allows to describe the formation of the SLE `bubbles' for $\kappa>4$. 
\end{abstract}

\section{Introduction}
The chordal Schramm-Loewner Evolution ($SLE(\kappa)$) is a one-parametric family of random curves in the upper half-plane. They were introduced by Schramm in order to describe the only possible conformally invariant scaling limits of interfaces in various lattice models from statistical physics. Their behaviour has been extensively studied and by now it is rather well understood in the stochastic sense. Many questions about their path-wise behaviour are still open. 

Given a curve, probably the first question that should be asked is whether this is a simple curve. In the case of SLE it was done by Rohde and Schramm \cite{schramm2005basic}. Among other things, they proved that SLE curves are getting rougher when $\kappa$ is increasing. Moreover, they have two phase transitions. For $\kappa\in[0,4]$ they are simple curves, for $\kappa\in(4,8)$ they have double points and form `bubbles' that cover the entire space, for $\kappa\ge 8$ they are space-filling curves. 

The standard proof of these phase transitions is based on the analysis of the so-called forward Loewner differential equation on the boundary of the half-plane. On the boundary it becomes a Bessel SDE with dimension $d>1$ and the first phase transition is equivalent to the phase transition of the Bessel SDE at $d=2$. For smaller dimension it hits the origin with probability one, for large dimensions it does not. 

In the present paper we base our analysis on the study of the \emph{backward} Loewner equation and its boundary behaviour which corresponds to the Bessel SDE with dimension below $1$. We express the question about the behaviour of SLE curves to the question whether the origin is an absorbing or reflecting boundary point for the Bessel SDE and whether this SDE has a unique solution when started from the origin. As a by-product of our analysis we can also express SLE `bubbles' in terms of the excursion decomposition of the Bessel process. 

We note that in the recent work \cite{shekhar2020complex} the authors study the existence of solutions of the complex Bessel SDE obtained from the backward Loewner differential equation driven by $\sqrt{\kappa}B_t$ for $\kappa \leq 4$. In contrast, in this work we study the phase transition in terms of uniqueness/non-uniqueness of solutions of this equation. Moreover, we obtain new structural information about the SLE curves for $\kappa>4$.

\section{Background}

Forward SLE maps $g_t$ are defined as solutions to the following ODE

\begin{equation}\label{5}
\partial_{t}g(t,z)=\frac{2}{g(t,z)-\sqrt{\kappa}B_{t}}\,, \hspace{10mm} g(0,z)=z, z \in \mathbb{H}\,
\end{equation}
where $B_t$ is the standard Brownian motion.  It is a well known fact  \cite{schramm2005basic} that with probability one $g_t$ is a continuous family of conformal maps from $\H_t$ to $\H$ where $\H_t$ is the unbounded component of the complement of a continuous curve $\gamma([0,t])$. This curve is called the SLE trace. Also, it is known \cite{schramm2005basic} that with probability one $g_t$ is continuous up to the boundary and that 
\[
\gamma(t)=\lim_{y\to 0} g_t^{-1}(\sqrt{\kappa}B_t+iy).
\]

Quite often it is more beneficial to work with inverse maps. Unfortunately, they satisfy a PDE which is much harder to analyse. Instead one can consider the so called backward SLE
\begin{equation}\label{6.1}
\partial_{t}h(t,z)=\frac{-2}{h(t,z)-\sqrt{\kappa}B_{t}}\,, \hspace{10mm} h(0,z)=z, z \in \mathbb{H}.
\end{equation}
For each fixed $t$, the map $h_t+\sqrt{\kappa}B_t$ has the same distribution as the map $g_t^{-1}$ (see Lemma $5.5$ of \cite{kemppainen2017schramm}). Moreover, \cite{rohde2016backward} it is possible to define the backward trace which, for a fixed time $t$ and up to a shift and reversal of parametrization, has the same distribution as the forward trace.

\subsection{Bessel processes}
Our results are based on the analysis of the boundary behaviour of $h_t$. On the real line the backward Loewner evolution \eqref{6.1} becomes a real SDE
\[
\partial_t h_t(x)=-\frac{2}{h_t(x)-\sqrt{\kappa}B_t}.
\]
Changing variables to $X_t=(h_t(x)-\sqrt{\kappa}B_2)/\sqrt{\kappa}$ we get
\begin{equation}
\label{eq:Bessel SDE}
d X_t=\frac{-2/\kappa}{X_t}dt+dB_t=\frac{d-1}{2}\frac{1}{X_t}dt+dB_t
\end{equation}
which is the Bessel SDE of dimension $d=d(\kappa)=1-4/\kappa$. It is easy to see that  the dimension $d \leq 0$ for $\kappa \leq 4$ and $d>0$ for $\kappa>4$.

It is known   that for  $\kappa \leq 4$ the origin is an absorbing point for the Bessel process and for $\kappa>4$ it is reflecting. In particular, in this case there is a unique non-negative instantaneously reflecting solution. Namely, \cite{aryasova2012properties} proves that there exists a unique strong non-negative solution to the Bessel equation that starts at zero and spends zero Lebesgue measure time there. We will use this solution to construct another strong solution that is not positive.

\section{Main result}

The first systematic study of SLE curves and their properties was performed by Rhode and Schramm \cite{schramm2005basic}. In particular, they proved that SLE curves are simple for $\kappa\in[0,4]$ and not simple for $\kappa>4$. Our main result is an alternative prove of this phase transition.  The first result is a new proof of a standard fact that SLE curves are simple when $\kappa\le 4$. 

The second result is more involved, beyond giving a new proof of the fact that the trace is not simple if $\kappa>4$ we describe how double points and the corresponding loops (or `bubbles') are formed.

%\begin{theorem}
%\label{thm: main theorem k>4}
%ADD PRECISE STATEMENT
%\end{theorem}

Before stating our main result, we introduce the following concepts.
\begin{definition}
For $m>0$, an $m$-macroscopic excursion is an excursion of the Bessel process of length at least $m >0$, where the length of an excursion is understood as the length of the time interval corresponding to the excursion.
%\begin{align}
%&dX_t=\frac{-2/\kappa}{X_t}dt+dB_t\,,\nonumber\\
%\end{align}
%started  from $X_0=0$, 
\end{definition}
\noindent Next, we introduce the notion of macroscopic hull.
\begin{definition}

A macroscopic hull is a compact set in $\bar{\mathbb{H}}$  with a simply connected complement, such that its intersection with the real axis is an interval $I \subset \mathbb{R}$ with strictly positive Lebesgue measure. We call the intersection of macroscopic hull with the real axis the base of the hull.
\end{definition}
\begin{definition}
Let $m>0$ be a real number. Let $\gamma(t)$ be the $SLE_{\kappa}$ curve. If there exists  two times $t_1$, $t_2 $ in the time interval such that $\gamma(t_1)=\gamma(t_2)$ and $|t_1-t_2| \geq m$, then we say that there is a self-touching of the curve after time at least $m$.
\end{definition}
We also introduce the following definition.
\begin{definition}
Let $m>0$. For the (backward) $SLE_{\kappa}$ trace,  an m-macroscopic double point is a double point that corresponds to self-touching of the trace after time at least $m>0$.
\end{definition}

The origin in $\mathbb{H}$ is a singularity for the backward Loewner differential equation, since the vector field $\frac{-2/\kappa}{Z_0}$ has infinite modulus. For almost every Brownian path, a notion of solution for $ t \in [0,T]$ for this differential equation is defined by the limit $Z_t(0):=\lim_{y \to 0+}Z_t(iy)$. 

%The existence of the limit and the continuity in time, is a consequence of the Rohde-Schramm Theorem (Theorem $2.2.2$ in the Introduction). We give a detailed explanation for this fact, in the following section of the chapter.%

We are now ready to state our main result.

%let $\phi(n)$ be a subpower function. Let us introduce the following 'boxes' in the upper-halfplane $\mathbb{H} :$
%$$A_{n,c, \phi}= \bigg\{ x+iy \in \mathbb{H}: |x| \leq \frac{\phi(n)}{\sqrt{n}}, \frac{1}{\sqrt{n}\phi(n)} \leq y \leq \frac{c}{\sqrt{n}} \bigg\}\,.$$ Then, for all fixed $t \in [0,1]\,,$ there is a unique solution for the backward Loewner differential equation started from the origin in $\mathbb{H}$, for almost every Brownian path. For defining the solution, let us consider a sequence of points $(z_n)_{n \in \mathbb{N}} \in A_{n,c,\phi}\,,$ such that $z_n \to 0$ as $n \to \infty\,.$ Then, for all fixed $t>0$, a way to define the solution is 
% $$ h_t(0+) := \lim_{n \to \infty}h_t(z_n) \,.$$ 

\begin{theorem}\label{result2}
For $\kappa\in (0,4]$ and for all times in a compact time interval, there is a unique (complex) solution of the backward Loewner differential equation started from the origin, a.s.. Moreover, for $\kappa \leq 4$, the backward $SLE_{\kappa}$ trace is a simple curve a.s.

For $\kappa >4$, there are at least two real solutions for the backward Loewner differential equation started from the origin, a.s.

Let $m>0$ be a positive real number. For $\kappa \in (4, \infty)$, on m-macroscopic excursions from the origin of the Squared Bessel process obtained from extensions of the backward $SLE_{\kappa}$ maps on the real line, we obtain macroscopic hulls with base depending on $m>0$ and $m$-macroscopic double points of the backward $SLE_{\kappa}$ trace. 
\end{theorem}
\begin{remark}
Compared with the proof of the a.s.\ simpleness of the $SLE_{\kappa}$ trace for $\kappa \leq 4$ that already exists in the literature, we show that this is implied by the behaviour at the origin the backward Bessel process obtained from the extension of the conformal maps to the boundary. 
\end{remark}

In addition, this analysis also offers a possible answer to the question: what pathwise properties of the Brownian motion influence the behavior of the $SLE_{\kappa}$ trace.  We will study this in more detail in the last section of this paper.

In a nutshell, in our analysis, we extend the conformal maps to the boundary and obtain that the dynamics of boundary points satisfy the Bessel SDE \eqref{eq:Bessel SDE}. We use known results about the solutions of the Bessel SDE  along with theorems about the backward Loewner evolution  in order to study   the backward Loewner evolution  started from the origin. 
For this, we use the result from \cite{aryasova2012properties}, in which there is a proof of the existence and uniqueness of strong positive solutions to the Bessel SDE started from the origin for dimension $ 0<d<1$. In our case, this gives the existence and uniqueness of a strong non-negative  solution for the Bessel SDE  on the real line for $\kappa>4$, since $d(\kappa)=1-\frac{4}{\kappa}.$ 
 
\begin{remark}
For dimensions $d \in (0,1]$ ($\kappa>4$) there exists a construction of the solutions to the SDE using Excursion Theory. Using this and the previous analysis, we obtain that beginning of m-macroscopic excursions for the Bessel SDE  give macroscopic hulls for the backward Loewner differential equation started from the origin. Thus, we obtain another structural information about the dynamics of the backward $SLE_{\kappa}$ hulls, and implicitly backward $SLE_{\kappa}$ traces. That is, the m-macroscopic excursions of the Bessel process started from the origin, create in this context, macroscopic hulls of the backward $SLE_{\kappa}$. In particular, when starting the backward Loewner differential equation from the beginning of m-macroscopic excursions of the Bessel processes, we obtain a macroscopic hull, a.s.
 For further details about the Excursion Theory construction of the Bessel processes for $d \in (0,1]$, see \cite{bertoin1990excursions}.
\end{remark}

\subsection{Backward Loewner differential equation and the Backward Loewner trace}

In the following we show how using the a.s. existence of the forward chordal $SLE_{\kappa}$,  we can obtain the a.s. existence of the backward Loewner trace.  For $\kappa>0$, the backward chordal $SLE_{\kappa}$ is defined by solving the backward Loewner differential equation with the driver $\sqrt{\kappa}B_t$, $0\leq t<\infty$. For any $t_0>0$, we have that the process
$(\sqrt{\kappa}B_{t_0-t}-\sqrt{\kappa}B_{t_0}, 0 \leq t\leq t_0) $ has the same distribution as $( \sqrt{\kappa}B_t$, $0\leq t \leq t_0)$. Rohde and Schramm showed in \cite{schramm2005basic} that for every $t_0 \in [0,T)$, $\sqrt{\kappa}B_{(t_0-t)}$, $0 \leq t \leq t_0$ generates a forward Loewner trace which we denote by $\beta_{t_0}(t_0-t), 0 \leq t \leq t_0$.  Then, using the identity in distribution for the Brownian motion from above, it is shown in \cite{rohde2016backward} that $\sqrt{\kappa}B_t$ generates the backward $SLE_{\kappa}$ traces $\beta_{t_0}$ for $0 \leq t_0\leq T.$

Note that $\beta_t$ is a continuous function defined on $[0,t]$. The parametrizations of the backward traces $\beta_t$ is different from the usual parametrization of the chordal $SLE_{\kappa}$ trace in the sense that
the backward traces $\beta_t$ is a continuous function defined on $[0,t]$ such that $\beta_0$ is the tip and $\beta_t$ is the root that is  an element of $\mathbb{R}\,.$ The difference with the chordal $SLE_{\kappa}$ is that the chordal trace is parametrized such that the root of it is $\beta_0\,.$
For further details about this construction, we refer to \cite{rohde2016backward}.

We also use in our analysis that for fixed time $T>0$, the law of the curves $\beta_T=\beta(0,T]$ generated by the backward Loewner differential equation is the same as the forward chordal $SLE_{\kappa}$ trace $\gamma(0,T]$ (modulo a translation with the driver $\sqrt{\kappa}B_T$).

\section{Proof of the main result}

\subsection{Proof of the Main Theorem for $\kappa\le 4$}

%Existence of the solution%
First, by an application of the Dominated Convergence Theorem one obtains that $\lim_{y \to 0+} h_t(iy)$ satisfies the backward Loewner differential equation driven by $\sqrt{\kappa}B_t$ for $\kappa \leq 4$ on compact time intervals, a.s.. 

Next, we show that SLE trace is simple for $\kappa\le 4$. One of the issues with the analysis of the Bessel process is that the drift coefficient blows up at the origin and many standard theorems about SDEs are not applicable. The standard way of dealing with this problem is to consider the square Bessel process $BESQ^\delta(x)$ which is defined as the unique strong solution of 
\[
X_t=x+2\int_0^t \sqrt{|X_s|}d B_s+\delta t.
\] 
It is well known that this process has the distribution of the square of the Bessel process of dimension $\delta$. 

For the square Bessel process, we have the following phase transition:
\begin{proposition}[Proposition $1.5$, Chapter XI of \cite{revuz2013continuous}]
For $\delta=0$, the point $0$ is absorbing and for $0<\delta<2$ the point $0$ is instantaneously reflecting.
\end{proposition}

Since the process with larger dimension stochastically dominates the process with a smaller dimension, this proposition immediately implies that the origin is absorbing for all $\delta<0$ as well.  This immediately implies that the origin is an absorbing point for the corresponding Bessel process.

%For $0<d<2$, the process is also reaching the origin, a.s. but the origin is not absorbing anymore.

To relate these results regarding the Bessel process with properties of the SLE trace we first recall a standard fact that SLE trace is simple with probability either $0$ or $1$ and that it is simple a.s. if and only if it does not hit the real line a.s. This follows from a simple fact that the trace $\gamma(t)$ is not simple if and only if there is $s$ such that $\gamma^s(t)$ hits the real line, where $\gamma^s$ is the trace of evolution driven by $\sqrt{\kappa}(B_{t+s}-B_s)$.
By scaling property of SLE it is enough to check whether SLE trace hits the real line by any fixed time, say time $1$. 

 In the next proof, we show that the last statement is equivalent to the statement that the origin is non-absorbing for the Bessel process.
 
To relate this fact to the Bessel process we need the following observation. Fix $s>0$, a deterministic time. We consider the Loewner chain $h_{s,s+t}$, i.e.  

\[
\partial_th_{s,s+t}(z)=\frac{-2}{ h_{s,s+t}(z)-\sqrt{\kappa}B_{s+t}}, \qquad h_{s,s}(z)=z.
\]

For any deterministic time shift $s>0$, when considering the continuous extension to the boundary, we obtain the same Bessel process with dimension $d=1-\frac{4}{\kappa}$ since the format of the shifted chain is the same with the one for $s=0$. This gives us that for any deterministic shift of time $s>0$, we obtain a Bessel process that has an absorbing origin. Thus, for a.e. Brownian Motion for all times $s \in [0,1]$ on these paths the corresponding Bessel process is absorbing. This excludes the possibility of random shifts of starting times $s \in [0,1]$ on which the Bessel process obtained from the continuous extension of the conformal maps to the real line started from those (random) times has any other solution than $X=0$, with positive probability.
\color{black}

For the backward Loewner differential equation, we have that on the real line the dynamics of points is governed by the following SDE

\begin{equation}\label{timerevert}
dZ_r=\frac{-2/\kappa}{Z_r}dr+d\tilde{B}_r\,,   \hspace{3mm} Z_0=x\in \mathbb{R}\setminus \{0\}.
\end{equation}
with dimensions $d \leq 0,$ where $\tilde{B_r}=B_1-B_{1-r}$.  
On the positive part of the real axis this is the equation of a Bessel process. In addition, in our analysis we have as a driver the same collection of Brownian paths, the dynamics on the negative part of the real axis is given by 
\begin{align}
d\tilde{X}_t&=\frac{-2}{\kappa\tilde{X}_t}dt+dB_t,\nonumber\\
\tilde{X}_0&=x_0<0\,.
\end{align}
Since $\tilde{X}_t$ is negative (until the first hitting time of zero), multiplying by $-1$ we obtain the following SDE
\begin{align}
dX_t&=\frac{-2}{\kappa X_t}dt-dB_t,\nonumber\\
X_0&=x_0>0\,,
\end{align}
with $X_t=-\tilde{X}_t\,.$ Thus, the dynamics on the negative part of the real axis, when coupling with the same Brownian drivers, is  equivalent to a dynamics on the positive side of the real axis with the driver $-B_t$ (that is still a Brownian motion).

If with positive probability the $SLE_{\kappa}$ trace were to be not simple, then  we obtain that with positive probability $\exists \tau(\omega)<1$, such that $\gamma^{(s)}(\tau) \cap\mathbb{R} \neq 0$.
Then, when considering the SDE \eqref{timerevert}, we obtain that the collection of Brownian Motion paths such that the origin is reflective for the equation \eqref{timerevert} has positive measure. Thus, when studying the maps $h_r(z)$ driven by $\sqrt{\kappa}B_{1-r},  r \in [0,1]$ we obtain a contradiction with the fact that $\{0\}$ is a.s. absorbing boundary : $Z_t \equiv 0$ for all $t >0$, for the backward Loewner differential equation for almost every Brownian path. Indeed if with positive probability the trace would touch the real line in $[0, 1]$, then we could find a positive mass of paths on which there exists times $r(\omega)$ such that $Z_{r(\omega)} \neq 0$, that gives a contradiction with the previous argument that shows that for almost every Brownian path for all times $t \in [0,1]$ on these paths the only possible solution for the Bessel process from the origin is $Z_t \equiv 0$. \\

This shows that the information on the boundary is enough to determine the behavior of the $SLE_{\kappa}$ trace. Indeed, since for $\kappa \in [0,4]$, for almost every Brownian path, the boundary behavior of the origin is absorbing, we obtain as a corollary that the the only possible notion of solution for the backward Loewner differential equation started from the origin is the one constructed before from 'inside the domain'.

\subsection{Proof of the Main Theorem for $\kappa >4 $}

In this section, we prove that in the case $\kappa >4$, almost surely there are at least two types of solutions for the backward Loewner differential equation started from the origin (i.e. solutions obtained from the same starting point and the same Brownian motion path).

In order to show this, we consider again the maps $h_t(z)$ solving the backward Loewner differential equation driven by the Brownian driver $\sqrt{\kappa}B_t$, for $ t \in [0,T]\,.$ We extend continuously the conformal mappings to the real line. 
We show that we can give meaning to a notion of solution to the Bessel SDE  that spends null Lebesgue measure time at the origin. Moreover, in the case $d>0$ (i.e. $\kappa>4$) we use this strong solution to construct a different solution, that gives the non-uniqueness in this regime.

%The main result that we use is proved in \cite{aryasova2012properties} and is phrased as follows: 

We use the notation  $L_a^{Z}(t)$ for the local time of the process $Z_t$ at the point $a \in \mathbb{R_+}$. Let us introduce the principal value correction for $d \in (0,1]$, i.e.
\begin{equation*}
Z_t = Z_0+ B(t) + \frac{d-1}{2}k(t)\,, Z_0 \geq 0,
\end{equation*}
where
\[
k(t) = P.V. \int_0^t\frac{1}{Z_s}ds :=\int_0^{\infty}a^{d-2}(L_a^{Z}(t)-L_{0}^{Z}(t))da.
\]
Aryasova and  Pilipenko proved in \cite{aryasova2012properties}   that the existence of a weak  solution spending zero time at the origin implies the existence and uniqueness of a non-negative strong solution spending zero time at the origin for the above equation. Note that this analysis allows us to view the strong solutions on the real line starting from the origin as functions of the Brownian motion paths.

We consider the following version of the previous equation (since the dynamics is defined on both the positive part of the real line and on the negative part of the real line).
We consider the following :
\begin{align}\label{SDEbig}
\tilde{k}(t)=P.V.\int_0^t\frac{dt}{Z_t}&=\int_{-\infty}^{\infty}sgn(a)|a|^{d-2}(L_a^Z(t)-L_0^Z(t))da.
\end{align}
Furthermore, we use the result from \cite{aryasova2012properties} to see this positive solution as a unique positive strong solution for the 'extended' SDE 
\begin{equation}\label{bigbigSDE}
Z_t = Z_0+ B(t) + \frac{d-1}{2}\tilde{k}(t)\,, Z_t \geq 0.
\end{equation}
Indeed this holds, since $Z_t \geq 0 $, gives $k(t)=\tilde{k}(t)$.
Note that this is a strong solution and thus can be thought as a measurable function $\phi: C(\mathbb{R_+})\to C(\mathbb{R_+})$ such that $\rho(t):=\phi(B_t(\omega))$ solves \eqref{SDEbig}. Note that this solution is unique, non-negative and starts from the origin, according to the result from \cite{aryasova2012properties}.

The next step is to consider $\tilde{B_t}(\omega)=-B_t(\omega)$, and as before $Z_0=0$.
We consider
\begin{align}
\tilde{\rho}(t):=-\phi(-B_t(\omega)).
\end{align}
We obtain that
\begin{align}
-\phi(-B_t(\omega))=\frac{d-1}{2}\tilde{k}(t)+B_t(\omega).
\end{align}
Thus, $\tilde{\rho}(t)$ is another solution to the equation (\ref{bigbigSDE}), with the same starting point and the same Brownian motion path. Note that the two solutions are different since one is positive and the other one is negative. Thus, for all times $ t \in [0,T]$, for $\kappa >4$, we have at least two different real solutions with the same starting point and with the same Brownian driver.
In particular, at each fixed time $t>0$, for $\kappa >4 $ we obtain for the backward Loewner differential equation started from the origin at least two solutions, a.s.

\section{Excursion Theory of the Bessel processes and behavior of the $SLE_{\kappa}$ traces}

\subsection{$SLE_{\kappa}$ traces and the Squared Bessel processes}
 
In this subsection, we discuss applications of the analysis developed so far on the $SLE_{\kappa}$ traces  for $\kappa >4 $.
We consider, as in \cite{lind2008holder}, the definition of the $SLE_{\kappa}$ trace using the mappings $(h_t(\sqrt{z})-\sqrt{\kappa}B_t)^2$, where $h_t(z):\mathbb{H} \to \mathbb H \setminus K_t$ are solving the backward Loewner differential equation.

The family of maps $(h_t(\sqrt{z})-\sqrt{\kappa}B_t)^2$ 
are conformal maps from $\mathbb C \setminus [0, \infty)$ to $\mathbb C \setminus [0, \infty)\setminus \tilde{K}_t$, where $\tilde{K}_t$ is the $SLE_{\kappa}$ hull corresponding to the maps  $$(h_t(\sqrt{z})-\sqrt{\kappa}B_t)^2.$$
Using the fact that the hulls of the backward SLE are locally connected, we obtain that the maps $(h_t(\sqrt{z})-\sqrt{\kappa}B_t)^2$ can be extended continuously to the real line.

When studying the extensions of the conformal maps satisfying the backward Loewner differential equation, we obtain on the real line two SDEs (that are both running on $\mathbb{R}_+$). These SDEs are obtained by applying It\^o's formula for $x \to x^2$. 
 
  These SDEs are the one satisfied by Squared Bessel processes.
Following \cite{revuz2013continuous}, by using the Yamada-Watanabe Theorem (see \cite{revuz2013continuous}), we obtain that the Squared Bessel process of any dimension, starting from $x\geq 0$ has a strong unique solutions a.s. 
 
Thus, in our case, we have the following SDEs. We call the first one, the upper Squared Bessel SDE and the second one the lower Squared Bessel SDE:
\begin{equation}\label{squaredBessel1}
d\tilde{Z_t}= \left(1-\frac{4}{\kappa}\right)dt+2\sqrt{\tilde{Z_t}}dB_t
\end{equation}

and 
\begin{equation}\label{squaredBessel2}
d\bar{Z_t}= \left(1-\frac{4}{\kappa}\right)dt-2\sqrt{\bar{Z_t}}dB_t.
\end{equation}

We couple these processes with the same Brownian driver $B_t(\omega)$, in the sense that we drive (\ref{squaredBessel1}) and (\ref{squaredBessel2}) with $B_t(\omega)$.
First, we restrict our attention to the solution of (\ref{squaredBessel1}). The analysis for the lower Squared Bessel process (i.e. the solution to $d\tilde{Z_t}= \left(1-\frac{4}{\kappa}\right)dt-2\sqrt{\tilde{Z_t}}dB_t$) will be performed in the same manner.

As in \cite{lind2008holder}, we denote the $SLE_{\kappa}$ trace in this new setting by $$\gamma_2(t) : =(g_t^{-1}(\sqrt{\kappa}B_t))^2.$$

%\red{Note that, by conformal invariance $\gamma_2(t)$ has the same law as the backward $SLE_{\kappa}$ trace in $\mathbb{H}$ (that modulo a shift has the same law with the forward $SLE_{\kappa}$ trace in $\mathbb{H}$). \qq{This is not true! Why do you need this statement?}}
%% We note that by definition, in this manner, the trace $\gamma_2(t)$ is parametrized by capacity.

For $d(\kappa)=1-\frac{4}{\kappa}>0$, there is a unique positive strong solution for this process starting from the origin such that $\tilde{Z}_t \geq 0$, a.s.. Also, for $2>d(\kappa)=1-\frac{4}{\kappa}>0$ the Squared Bessel process is recurrent. Thus, we can apply elements of Excursion Theory.
% developed in the previous section.

\subsection{The m-macroscopic excursions of the real Squared  Bessel process and application to the study of the backward $SLE_{\kappa}$ traces}

Fix $m>0$. In this section, we analyze the backward Loewner flow during a m-macroscopic excursion from the origin of the (upper) Squared Bessel process on $\mathbb{R}$. The same analysis holds for the (lower) Squared Bessel process.

Let us consider the collection of m-macroscopic excursions of the solution to the real Squared Bessel SDE. Let us consider the backward Loewner differential equation from the origin, for time $t \in [0,\infty).$
For the Squared Bessel process which is obtained by extending the conformal maps on the real line, we have the following result (see Section $4.6$ of \cite{rogers2000diffusions}), for $0<d<2$ 
$$\int_0^{\infty}\bold{1}_{\{\tilde{Z}_t=0\}}dt=0.$$ 
Thus, the origin is instantaneously reflecting for this process for all the times $[0,\infty)$.
Among these times, there exist a set of times which correspond to beginning of m-macroscopic excursions of the Squared Bessel process from the origin.
 In order to see this, we prove the following result:
 $$\mathbb{P}(\text{there is an excursion of length at least} \hspace{2mm} m)=1.$$

Since the Bessel process for $d >0$ is not identically zero, with positive probability there is at least one excursion interval with length $l>m$. The choice of this constant $m>0$ is arbitrary and the modification of this constant does not change the result.  Following the approach in \cite{bertoin1998levy}, let us consider the event $$\Lambda_t=\{ \text{all the excursion intervals with right-end point d} <t \hspace{2mm} \text{ have length} \hspace{2mm} l \leq m \}.$$
 If the first return to $\{0\}$ of the process (that is a stopping time) is infinite, then the Markov process has an excursion interval of length $l>m$.  If this stopping time is finite, then we apply the Markov property and we get that $\mathbb{P}(\Lambda_{3t})\leq \mathbb{P}(\Lambda_{t})^2$ and by iteration we get that $\mathbb{P}(\Lambda_{3^nt})\leq \mathbb{P}(\Lambda_{t})^{2n}$ 
for $n \in \mathbb N$. Thus, we get that $\lim_{s\to \infty}\mathbb{P}(\Lambda_s)=0$. 
 
We further investigate the continuous extensions $(h_t(\sqrt{z})-\sqrt{\kappa}B_t)^2$, for times $t \in [0, \infty)$. Among these times, there will be times $r=r(\omega)$ (that depend on the path of the Brownian motion) that are beginning of m-macroscopic excursions of the real Squared Bessel process, started from the origin. Let us consider 
\[
\tilde{h}_r(z):=h_{s+r} \circ h_r^{-1}(z),
\]
 for $s\geq 0$. In particular,  we obtain that the family of extended conformal maps $(\tilde{h}_r(0)-\sqrt{\kappa}\hat{B}_r)^2$, where $\hat{B}_r:=B_{t+r}-B_r$ for times in $ [r(\omega), r(\omega)+m)$, map the origin to the real line. In order to see this, let us assume that the image of the origin under the maps $(\tilde{h}_r(0)-\sqrt{\kappa}\hat{B}_r)^2$ is  in $\mathbb H$. Then, for any $\tilde{m} \in (0, m)$ we would have two images of the origin, for the extended conformal maps, one in $\mathbb{H}$ and the unique strong solution of the Squared Bessel SDE at time $r(\omega) + \tilde{m}$.

The maps $\tilde{h}_r(z)$ are continuously extended to the boundary (since they are obtained from compositions of maps that are continuously extended to the boundary). Hence, since the maps $\tilde{h}_r(z)$ should map the origin to a unique point, we obtain a contradiction. Thus, the two images of the origin should coincide. Since the backward $SLE_{\kappa}$ trace grows from the origin, we obtain that when starting the backward Loewner flow from the set of m-macroscopic excursions of the underlying Bessel, we obtain points of intersection of the trace with the real line. Using that a.s. there is an excursion interval of length at least $m>0$ as shown before, we obtain that the closure of backward $SLE_{\kappa}$ hulls are macroscopic hulls whenever we have m-macroscopic excursions of the Squared Bessel process. Moreover, since the Squared Bessel process is recurrent and the backward $SLE_{\kappa}$ trace starts from the origin, we obtain also double points (that correspond to self-touchings of the curve after at least $m>0$ time).

Moreover, from topological considerations, we obtain at the end of a m-macroscopic excursion of the unique strong solution of the Squared Bessel process, a closed `large bubble' of the backward SLE trace formed from the outer boundary of a portion of the curve, where $t (\omega)$ is the beginning of a large excursion that collects the self intersections of the backward $SLE_{\kappa}$ trace with itself and with real line, that happened on $[t(\omega),t(\omega)+ m)$. 
In contrast with the forward flow approach, this approach not only recovers the phase transition in the behaviour of the $SLE_{\kappa}$ traces, but also it provides structural information about the backward $SLE_{\kappa}$ traces.

The same analysis holds also for the solution to the SDE (\ref{squaredBessel2}). When driving them simultaneously with $B_t(\omega)$ and $-B_t(\omega)$  we obtain the complete picture of the evolution of the backward Loewner hulls in $\mathbb{C}\setminus [0, \infty)$.

\subsection{Brownian motion and the behavior of the $SLE_{\kappa}$ trace for $\kappa \in (4,\infty)$}
 
In addition, this analysis also offers a possible answer to the question: what pathwise properties of the Brownian motion influence the behavior of the backward $SLE_{\kappa}$ trace? The previous analysis suggests that the classes of points along a Brownian driver that allow the process to escape the origin (and to further reflect at the origin), give the possibility of the non-simpleness of the backward $SLE_{\kappa}$ traces for $\kappa>4$. Compared with the case when the parameter $\kappa \in [0,4]$, where the solution along the positive real line is stuck at the origin, in the case $\kappa \in (4, \infty)$ there is a notion of strong solution along the real line. Equivalently, when the parameter $\kappa \in [0,4]$ is not large enough to allow the solution of the backward Bessel SDE to escape the origin, there is no solution along the real line. 

In general, the change in behavior of the backward Bessel process started from the origin, gives also the following additional information: the collection of points that represent the beginning of macroscopic excursions of the real Bessel process of dimensions $d \in (0,1)$ that is obtained from the extensions of the conformal maps to the real line, give also the formation of macroscopic bubbles for the backward $SLE_{\kappa}$ hulls.

\bibliographystyle{plain}
\bibliography{literature}
%\nocite*

\end{document}